\documentclass[11pt]{article} 

\usepackage[utf8]{inputenc} 
\usepackage{comment}

\usepackage[numbers]{natbib}
\usepackage{hyperref}
\usepackage{geometry} 
\usepackage{graphicx} 
\usepackage{booktabs}
\usepackage{array} 
\usepackage{paralist} 
\usepackage{verbatim}
\usepackage{bm} 
\usepackage{subfig} 
\usepackage{amsmath}

\usepackage{fancyhdr} 
 
\usepackage{sectsty}
\allsectionsfont{\sffamily\mdseries\upshape} 
\usepackage[nottoc,notlof,notlot]{tocbibind} \usepackage[titles,subfigure]{tocloft}

\numberwithin{equation}{section}

\title{Derivation of the stochastic Hamilton-Jacobi-Bellman equation}
\author{
    Vasil Yordanov\\
    Sofia University St. Kliment Ohridski\\
    \texttt{v.yordanov@phys.uni-sofia.bg}
}

\begin{document}
\maketitle
\abstract 
In the present paper, we provide a detailed derivation of the stochastic Hamilton–Jacobi–Bellman (HJB) equation.

\section{Introduction}
In this paper, we will derive the stochastic Hamilton–Jacobi–Bellman (HJB) equation. The derivation presented here is inspired by the following papers:~\cite{Kappen2005,Kappen2011,Fleming2006}.

\section{Derivation of the stochastic HJB equation}
The stochastic equation of motion of the particle is:
\begin{equation}
\label{eq:StochasticProcess_appendix}
dx_\mu=u_\mu ds+ \sigma_\mu dW_\mu, \quad \mu=0..3,
\end{equation}
where $x_\mu$ are the spacetime coordinates $\mathbf  x$ of the particle, and $u_\mu$ are the components of the four-velocity $\mathbf u$. 

The action is postulated as minimum of the expected value of stochastic action:
\begin{equation}
\label{eq:eq_for_cost}
S(\mathbf x_i, \mathbf u(\tau_i \rightarrow \tau_f))= \min_{\mathbf u(\tau_i \rightarrow \tau_f)} \left< \int_{\tau_i}^{\tau_f} ds\, \mathcal{L}(\mathbf x(s), \mathbf u(s),s) \right>_{\mathbf  x_i},
\end{equation}
where $\mathcal{L}(\mathbf x(s), \mathbf u(s), s)$ is the Lagrangian of the test particle, which is a function of the control policy $\mathbf u(s)$ and four-coordinates $\mathbf x(s)$ at proper time $s$ .  
The subscript $\mathbf x_i$ on the expectation value means that the expectation is
over all stochastic trajectories that start at $\mathbf x_i$.

The task of optimal control theory~\cite{bellman1954} is to find the control $\mathbf u(s)$, $\tau_i < s < \tau_f$, denoted as $\mathbf u(\tau_i \rightarrow \tau_f)$, that minimizes the expected value of the action $S(\mathbf x_i, \mathbf u(\tau_i \rightarrow \tau_f))$.

We introduce the optimal cost-to-go function for any intermediate proper time $\tau$, where  $\tau_i < \tau < \tau_f$:
\begin{equation}
\label{eq:cost_to_go_function}
J(\tau,  \mathbf  x_\tau)=\min_{\mathbf u(\tau  \rightarrow \tau_f )} \left< \int_{\tau}^{\tau_f} ds\, \mathcal{L}(s, \mathbf x_s, \mathbf u_s)  \right>_{ \mathbf x_\tau}
\end{equation}

By definition, the action $S(\mathbf x_i, \mathbf u(\tau_i \rightarrow \tau_f))$ is equal to the cost-to-go function $J(\tau_i,  \mathbf  x_{\tau_i})$ at the initial proper time and spacetime coordinate:
\begin{equation}
S(\mathbf x_i, \mathbf u(\tau_i \rightarrow \tau_f))=J(\tau_i,  \mathbf  x_{\tau_i})
\end{equation}

We can rewrite recursive formula for $J(\tau, \mathbf x_\tau)$ for any intermediate time $\tau'$, where $\tau < \tau' < \tau_f$:
\begin{equation}
\begin{aligned}
J(\tau,  \mathbf x_\tau)
&=\min_{\mathbf u(\tau  \rightarrow \tau' )} \left< \int_\tau^{\tau'} ds\,\mathcal{L}(s, \mathbf x_s, \mathbf u_s) + \int_{\tau'}^{\tau_f} ds\,\mathcal{L}(s, \mathbf x_s, \mathbf u_s) \right>_{\mathbf x_\tau} \\
&=\min_{\mathbf u(\tau  \rightarrow \tau' )} \left< \int_\tau^{\tau'} ds\,\mathcal{L}(s, \mathbf x_s, \mathbf u_s) + \min_{\mathbf u(\tau'  \rightarrow \tau_f )} \left< \int_{\tau'}^{\tau_f} ds\,\mathcal{L}(s, \mathbf x_s, \mathbf u_s)  \right>_{\mathbf x_{\tau'}} \right>_{\mathbf x_\tau} \\
&=\min_{\mathbf u(\tau  \rightarrow \tau' )} \left< \int_\tau^{\tau'} ds\,\mathcal{L}(s, \mathbf x_s, \mathbf u_s) + J(\tau',  \mathbf x_{\tau'}) \right>_{\mathbf x_\tau}.
\end{aligned}
\end{equation}

In above equation we split the minimization over two intervals. These are not independent, because the second minimization is conditioned on the starting value $x_{\tau'}$, which depends on the outcome of the first minimization.

If $\tau'$ is a small increment of $\tau$, $\tau'=\tau + d\tau$ then:
\begin{equation}
\label{eq:J_recursivly_expressed_by_L}
J(\tau,  x_\tau)=\min_{\mathbf u(\tau  \rightarrow \tau+d\tau )} \left<\mathcal{L}(\tau, \mathbf x_\tau, \mathbf u_\tau) d \tau + J(\tau+d\tau, \mathbf x_{\tau + d\tau}) \right>_{\mathbf x_\tau}
\end{equation}

We must take a Taylor expansion of $J$ in $d\tau$ and $d \mathbf x$. However, since $\left<d \mathbf x^2 \right>=\sigma^2 d\tau$ is of order $d\tau$, we must expand up to order $d \mathbf x^2$:

\begin{equation}
\label{eq:avg_cost-to-go}
\begin{aligned}
& \left< J(\tau+d\tau,   \mathbf x_{\tau + d\tau}) \right>_{\mathbf x_\tau} 
=\int d  \mathbf x_{\tau+d\tau} \mathcal N (\mathbf x_{\tau + d\tau} |  \mathbf x_\tau, \sigma d\tau) J(\tau + d\tau,   \mathbf x_{\tau+d\tau}) \\
&=\int d  \mathbf x_{t+d\tau} \mathcal N (\mathbf x_{t+d\tau} |  \mathbf x_\tau, \sigma d\tau) (J(\tau,  \mathbf x)+d\tau \partial_\tau J(\tau,  \mathbf x_\tau) +dx^\mu \partial_\mu J(\tau,  \mathbf x_\tau) + dx^\mu dx^\nu \frac{1}{2} \partial_{\mu \nu} J(\tau,  \mathbf x_\tau) ) \\
&=J(\tau, \mathbf x)+ d\tau \partial_\tau J(\tau, \mathbf x_\tau) + \int d  x_{t+d\tau} \mathcal N (\mathbf x_{t+d\tau} |  \mathbf x_\tau, \sigma d\tau) (dx^\mu \partial_\mu J(\tau,  \mathbf x_\tau) + dx^\mu dx^\nu \frac{1}{2} \partial_{\mu \nu} J(\tau,  \mathbf x_\tau) ) \\
&=J(\tau, \mathbf x)+d\tau \partial_\tau J(\tau, \mathbf x_\tau) + \left<dx^\mu \right> \partial_\mu J(\tau, \mathbf x_\tau) + 
\frac{1}{2} \left<dx^\nu dx^\mu \right> \partial_{\nu \mu} J(\tau, \mathbf x_\tau)
\end{aligned}
\end{equation}

Here $\mathcal N (\mathbf x_{\tau + d\tau}|\mathbf x_\tau, \sigma d\tau)$ is the conditional probability starting from state $\mathbf x_\tau$ to end up in state $\mathbf x_{\tau+d\tau}$. The integration is over the entire spacetime.
In the above equation we also use the notation $ \partial_\mu J(\tau, \mathbf x)$ as partial derivative with respect to $x_\mu$.

We can calculate $\left< dx^\mu \right>$ using equations~\eqref{eq:StochasticProcess_appendix}:

\begin{equation}
\begin{aligned}
&\left< dx^\mu \right> = \int d \mathbf x_{\tau + d\tau} \mathcal N (\mathbf x_{\tau + d\tau} |  \mathbf x_\tau, \sigma^\mu d\tau) dx^\mu \\
&= \int d  \mathbf x_{\tau + d\tau} \mathcal N (\mathbf x_{\tau + d\tau} | \mathbf x_\tau, \sigma d\tau) (\mathbf u^\mu d\tau + \sigma^\mu d W^\mu) \\
&=\mathbf u^\mu d\tau \int d  \mathbf x_{t+d\tau} \mathcal N (\mathbf x_{\tau + d\tau} |  \mathbf x_\tau, \sigma^\mu d\tau) +  \int d  \mathbf x_{\tau + d\tau} \mathcal N (\mathbf x_{\tau + d\tau} |  \mathbf x_\tau, \sigma^\mu d\tau) \sigma^\mu d W^\mu
\end{aligned}
\end{equation}

From where:
\begin{equation}
\left<dx^\mu \right> = u^\mu d\tau
\end{equation}

In similar way we calculate $\left<dx^\nu dx^\mu \right> $:
\begin{equation}
\begin{aligned}
& \left< dx^\nu dx^\mu \right>  = \int d  \mathbf x_{\tau+d\tau} \mathcal N (\mathbf x_{\tau + d\tau} |  \mathbf x_\tau, \sigma^\mu d\tau) dx^\mu dx^\nu \\
&= \int d  \mathbf x_{\tau + d\tau} \mathcal N (\mathbf x_{\tau + d\tau} |  \mathbf x_\tau, \sigma^\mu d\tau) (u^\mu d\tau + \sigma^\mu d W^\mu) (u^\nu d\tau + \sigma^\nu d W^\nu) \\ 
&= \int d  \mathbf x_{\tau + d\tau} \mathcal N (\mathbf x_{\tau + d\tau} |  \mathbf x_\tau, \sigma^\mu d\tau) (u^\mu u^\nu d^2 \tau + u^\mu \sigma^\nu d\tau d W^\nu + \sigma^\mu d W^\mu u^\nu d\tau+  \sigma^\mu d W^\mu  \sigma^\nu d W^\nu ) \\ 
\end{aligned}
\end{equation}

From where we derive:
\begin{equation}
\left< dx^\nu dx^\mu \right> = 0, \mu \neq \nu,\quad  \left< (dx^\mu)^2 \right> = \sigma^\mu \sigma^\mu d \tau
\end{equation}

After substituting the above equations in~\eqref{eq:avg_cost-to-go} we derive the stochastic HJB equation:
\begin{equation}
\label{eq:HJB}
-\partial_\tau J(\tau, \mathbf x)=\min_{\mathbf u} \left( \mathcal{L}(\tau, \mathbf x, \mathbf u) + u^\mu \partial_\mu J(\tau, \mathbf x) + \frac{1}{2} \sum_{\mu=0}^{3} \sigma^{\mu} \sigma^\mu \partial_{\mu \mu} J(\tau, \mathbf x_\tau) \right)
\end{equation}

It is clear the from its definition of $J(\tau, x_\tau)$ that we have the following boundary condition:
\begin{equation}
\label{eq:boundary_condition}
J(\tau_f, x_{\tau_f}) = 0
\end{equation}

The optimal control at the current $\mathbf x$, $\tau$, is given by:
\begin{equation}
u(\mathbf x, \tau) = \arg \min_{\mathbf u} \left(\mathcal{L}(\tau,  \mathbf x, \mathbf u) + u^\mu \partial_\mu J(\tau, \mathbf x) \right).
\end{equation}

\bibliographystyle{unsrt}
\bibliography{refs}

\begin{thebibliography}{1}

\bibitem{Kappen2005}
Hilbert~J. Kappen.
\newblock Path integrals and symmetry breaking for optimal control theory.
\newblock {\em Journal of Statistical Mechanics: Theory and Experiment},
  2005(11):P11011--P11011, nov 2005.

\bibitem{Kappen2011}
Hilbert Kappen.
\newblock Optimal control theory and the linear bellman equation.
\newblock {\em Neural Networks}, 08 2011.

\bibitem{Fleming2006}
Wendell~H. Fleming and H.M. Soner.
\newblock {\em Controlled Markov Processes and Viscosity Solutions}.
\newblock Stochastic Modelling and Applied Probability. Springer, New York, NY,
  2 edition, 2006.
\newblock Number of Pages: XVII, 429.

\bibitem{bellman1954}
Richard Bellman.
\newblock The theory of dynamic programming.
\newblock {\em Operations Research}, 2(3):275--285, 1954.

\end{thebibliography}

\end{document}